\documentclass[11pt,leqno]{smfart}

\usepackage[latin1]{inputenc}
\usepackage[english,francais]{babel}
\usepackage{amsmath} 
\usepackage{amssymb}
\usepackage{a4wide}
\usepackage{epsfig}
\usepackage{enumerate}
\usepackage{stmaryrd}
\usepackage[T1]{fontenc}
\usepackage[sans]{dsfont}

\newtheorem{theoreme}{Th\'eor\`eme}%[section]

\newtheorem{corollaire}[theoreme]{Corollaire}

\renewcommand{\theequation}{\arabic{equation}}
\newcounter{hypo}
\newcounter{hyphyp}

\def\C{{\mathbb C}}
 
\def\R{{\mathbb R}}

\def\CO{\mathcal {O}}

\def\CS{\mathcal {S}}

 \def\im{\mathop{\rm Im}\nolimits}

\def\Op{\mathop{\rm Op}\nolimits}

\def\arg{\mathop{\rm arg}\nolimits}

\def\one{\mathds{1}}
\def\<{\langle}
\def\>{\rangle}

\renewcommand{\theenumi}{\sl{\roman{enumi}}}
\renewcommand{\labelenumi}{\theenumi)}

\makeatletter
 \@addtoreset{equation}{section}
 \makeatother
 \renewcommand{\theequation}{\thesection.\arabic{equation}}

\begin{document}

\author[J.-F. Bony]{Jean-Fran\c{c}ois Bony}
\address[J.-F. Bony]{IMB, CNRS (UMR 5251), Universit\'e de Bordeaux, 33405 Talence, France}
\email[J.-F. Bony]{bony@math.u-bordeaux.fr}

\author[S. Fujiie]{Setsuro Fujii\'e}
\address[S. Fujiie]{Department of Mathematical Sciences, Ritsumeikan University, 1-1-1 Noji-Higashi, Kusatsu, 525-8577 Japan}
\email{fujiie@fc.ritsumei.ac.jp}

\author[T. Ramond]{Thierry Ramond}
\address[T. Ramond]{LMO (UMR CNRS 8628), Universit\'e Paris Sud 11, 91405 Orsay, France}
\email{thierry.ramond@math.u-psud.fr}

\author[M. Zerzeri]{Maher Zerzeri}
\address[M. Zerzeri]{Universit\'e Paris 13, Sorbonne Paris Cit\'e, LAGA, CNRS (UMR 7539), 93430 Villetaneuse, France}
\email{zerzeri@math.univ-paris13.fr}



\title{Propagation des singularit\'es et r\'esonances}

%\date{}

\begin{abstract}
Dans le cadre de l'\'etude des r\'esonances semiclassiques, on pr\'ecise le lien entre majoration polynomiale du prolongement de la r\'esolvante et propagation des singularit\'es \`a travers l'ensemble capt\'e. Cette approche permet d'\'eliminer l'infini et de concentrer l'\'etude pr\`es de l'ensemble capt\'e. Nous l'avons utilis\'ee dans des travaux ant\'erieurs pour obtenir l'asymptotique des r\'esonances dans diverses situations g\'eom\'etriques.
\end{abstract}

\begin{altabstract}
{\bf Propagation of singularities and resonances.}
In the framework of semiclassical resonances, we make more precise the link between polynomial estimates of the extension of the resolvent and propagation of the singularities through the trapped set. This approach makes it possible to eliminate infinity and to concentrate the study near the trapped set. It has allowed us in previous papers to obtain the asymptotic of resonances in various geometric situations.
\end{altabstract}

\maketitle

\section{Introduction}

Dans cette note, on consid\`ere un op\'erateur de Schr\"odinger semiclassique sur $L^{2} ( \R^{n} )$, $n \geq 1$,
\begin{equation}\label{a0}
P = -h^{2} \Delta + V ( x ) ,
\end{equation}
avec $V \in C^{\infty}_{0} ( \R^{n} ; \R )$. Plus g\'en\'eralement, on peut traiter le cas des op\'erateurs pseudodiff\'erentiels rentrant dans le cadre des op\'erateurs ``black box'' de Sj\"{o}strand et Zworski (voir par exemple \cite{Sj01_02}). La r\'esolvante $( P - z )^{- 1}$ est alors bien d\'efinie pour $\im z > 0$ et admet, en tant qu'op\'erateur de $L^{2}_{\rm comp} ( \R^{n} )$ dans $L^{2}_{\rm loc} ( \R^{n} )$, un prolongement m\'eromorphe  \`a un secteur angulaire $\{ z \in \C ; \ - \theta_{0} < \arg z \leq 0 \}$ avec $\theta_{0} > 0$. On appelle r\'esonances de $P$ les p\^oles de ce prolongement. Nous renvoyons le lecteur \`a \cite{Sj01_02} pour une pr\'esentation d\'etaill\'ee de cette th\'eorie.

On fixe une \'energie $E_{0} > 0$ et on s'int\'eresse aux r\'esonances $z=z_h$ dans $B ( E_{0} , C h \vert \ln h \vert ) = \{ z \in \C ; \ \vert z - E_{0} \vert < C h \vert \ln h \vert \}$ avec $C > 0$, qui est le domaine naturel d'\'etude des r\'esonances dans le cadre $C^{\infty}$ (voir par exemple \cite{Ma02_01}). On dit que {\it la r\'esolvante tronqu\'ee de $P$ v\'erifie une estimation polynomiale dans} $\Omega =\Omega_h\subset B ( E_{0} , C h \vert \ln h \vert)$ lorsque $P$ n'a pas de r\'esonance dans $\Omega$ et, pour tout $\chi \in C^{\infty}_{0} ( \R^{n} )$, il existe $N > 0$ tel que
\begin{equation} \label{a2}
\big\Vert \chi ( P - z )^{- 1} \chi \big\Vert \lesssim h^{- N} ,
\end{equation}
uniform\'ement pour $z \in \Omega$. Ici,  $( P - z )^{- 1}$ d\'esigne le prolongement de la r\'esolvante mentionn\'e ci-dessus. Pour que cela soit vrai, il suffit d'ailleurs que \eqref{a2} soit satisfaite pour une seule fonction $\chi \in C^{\infty}_{0} ( \R^{n} )$ \'egale \`a $1$ dans un compact assez grand. N\'eanmoins $N$ d\'epend de $\chi$. En fait, cette propri\'et\'e est \'equivalente \`a une majoration polynomiale de la r\'esolvante de $P_{\theta}$, l'op\'erateur distordu d'angle $\theta = M h \vert \ln h \vert$ avec $M \gg 1$. Ce genre d'estimation correspond donc \`a l'absence de pseudospectre semiclassique polynomial pour l'op\'erateur non-autoadjoint $P_{\theta}$. Les remarques pr\'ec\'edentes d\'ecoulent de la proposition D.1 de \cite{BoFuRaZe16_01} \'etendue au cas $\vert \im z \vert \lesssim h \vert \ln h \vert$.

Soit $p ( x , \xi ) = \xi^{2} + V ( x )$ le symbole semiclassique de $P$ et
$H_p=\partial_\xi p\cdot \partial_x - \partial_xp\cdot\partial_\xi =2\xi\cdot\partial_x -\nabla V\cdot\partial_\xi$
son champ hamiltonien. Pour $E \in \R$, on note
\begin{equation*}
\Gamma_{\pm} ( E ) = \{ \rho \in p^{- 1} ( E ) ; \ \exp ( t H_{p} ) ( \rho ) \not\to \infty \text{ quand } t \to \mp \infty \} ,
\end{equation*}
les ensembles sortant/entrant \`a l'infini.
L'ensemble capt\'e \`a \'energie $E$ est d\'efini par $K ( E ) = \Gamma_{-} ( E ) \cap \Gamma_{+} ( E )$. On sait que $\Gamma_{\pm} ( E )$ est ferm\'e et que $K ( E )$ est compact (voir \cite[Appendix]{GeSj87_01}).

Toutes les fonctions $u = u_{h} \in L^{2} ( \R^{n} )$ consid\'er\'ees dans la suite seront suppos\'ees polynomialement born\'ees, c'est \`a dire telles qu'il existe $N \in \R$ avec $\Vert u \Vert  \lesssim h^{- N}$. Pour $K \subset T^{*} \R^{n}$ un compact et $u$ une fonction polynomialement born\'ee, on dit que $u = 0$ microlocalement pr\`es de $K$ s'il existe $\varphi \in C^{\infty}_{0} ( T^{*} \R^{n} )$ avec $\varphi = 1$ pr\`es de $K$ tel que $\Op ( \varphi ) u = \CO ( h^{\infty} )$ en norme $L^{2}$. Le lecteur pourra se reporter \`a Dimassi et Sj{\"o}strand \cite{DiSj99_01}, Martinez \cite{Ma02_02} ou Zworski \cite{Zw12_01} pour une pr\'esentation d\'etaill\'ee de l'analyse microlocale semiclassique. On consid\`ere le {\it probl\`eme de Cauchy microlocal pr\`es de l'ensemble capt\'e}, c'est \`a dire les solutions $u$ de 
\begin{equation} \label{a3}
\left\{ \begin{aligned}
&( P - z ) u = v &&\text{ microlocalement pr\`es de } K ( E_{0} ) ,  \\
&u = u_{0} &&\text{ microlocalement pr\`es de tout point de } \Gamma_{-} ( E_{0} ) \setminus K ( E_{0} ) .
\end{aligned} \right.
\end{equation}
Cette \'equation doit \^etre vue comme une \'equation de propagation, $u_{0}$ \'etant la donn\'ee initiale et $v$ le second membre. On parle d'unicit\'e du probl\`eme de Cauchy microlocal si toute solution $u$ de \eqref{a3} avec $v = u_{0} = 0$ est microlocalement nulle pr\`es de $K ( E_{0} )$. Dans la litt\'erature, ce genre de r\'esultat porte parfois le nom de propagation des singularit\'es. On parle d'existence du probl\`eme de Cauchy microlocal si, pour toutes donn\'ees $v , u_{0}$ v\'erifiant $( P - z ) u_{0} = v$ microlocalement pr\`es de tout point de $\Gamma_{-} ( E_{0} ) \setminus K ( E_{0} )$, il existe une fonction $u$ solution de \eqref{a3}.

Les notions introduites ci-dessus sont reli\'ees par le r\'esultat suivant.

\begin{theoreme}\sl \label{a1}
Sous les hypoth\`eses pr\'ec\'edentes, sont \'equivalentes:
\renewcommand{\theenumi}{\roman{enumi}}
\begin{enumerate}
\item une estimation polynomiale dans $\Omega$ de la r\'esolvante tronqu\' ee, \label{a4}
\item l'unicit\'e du probl\`eme de Cauchy microlocal \eqref{a3} pour $z \in \Omega$. \label{a5}
\end{enumerate}
\end{theoreme}

L'implication \eqref{a5}$\Rightarrow$\eqref{a4} a d\'ej\`a \'et\'e prouv\'ee dans \cite{BoFuRaZe16_01}. Elle s'est r\'ev\'el\'ee un argument cl\'e pour le calcul asymptotique des r\'esonances dans diverses situations g\'eom\'etriques (cf. \cite{BoFuRaZe16_02,BoFuRaZe16_01}). La r\'eciproque est d\'emontr\'ee dans la section suivante. L'id\'ee est d'\'etendre une solution du probl\`eme de Cauchy microlocal pr\`es de l'ensemble capt\'e, \`a l'aide du flot quantique associ\'e \`a l'op\'erateur distordu $P_{\theta}$, afin de construire un quasimode de $P_{\theta} - z$.

A notre connaissance, dans tous les cas o\`u la distribution des r\'esonances a \'et\'e \'etablie, la r\'esolvante tronqu\'ee est polynomialement born\'ee d\`es que l'on s'\'eloigne des r\'esonances. Par ailleurs, Stefanov a d\'emontr\'e dans \cite{St01_01} sous des hypoth\`eses tr\`es g\'en\'erales que
\begin{equation*}
\int_{[E_{0} - \delta , E_{0} + \delta]} \big\Vert \chi ( P - z )^{- 1} \chi \big\Vert \, d z \lesssim h^{- N} ,
\end{equation*}
avec $N , \delta > 0$ et donc \eqref{a4} est presque toujours vrai sur l'axe r\'eel. Mais \eqref{a4} est faux \`a distance $\CO ( h^{\infty} )$ des r\'esonances. Ainsi, les assertions \'equivalentes du th\'eor\`eme peuvent \^etre vraies ou fausses pour $\vert \im z \vert \lesssim h \vert \ln h \vert$. Par contre, il n'est pas possible d'avoir une estimation polynomiale de la r\'esolvante tronqu\'ee plus profond\'ement dans le complexe (i.e. $\im z \ll - h \vert \ln h \vert$), comme le montre la proposition 1.5 de \cite{BoPe13_01}.

Pour prouver l'existence de r\'esonances tr\`es proches de l'axe r\'eel, il est souvent fait usage de quasimodes de $P - z$ \`a support compact, comme par exemple dans Tang et Zworski \cite{TaZw98_01}. De tels quasimodes n'existent pas lorsque $\im z$ est d'ordre $h$ (voir \cite{BoPe13_01}). Le th\' eor\` eme \ref{a1} permet de s'en passer, en transformant tout d\'efaut d'unicit\'e du probl\`eme de Cauchy microlocal en la production d'un quasimode global. De plus, comme le montre le r\'esultat suivant, les deux implications du th\'eor\`eme \ref{a1} permettent parfois d'obtenir l'existence de r\'esonances.

\begin{corollaire}\sl \label{a16}
Soit $D \subset B ( E_{0} , C h \vert \ln h \vert)$ un disque ouvert tel que l'unicit\'e du probl\`eme de Cauchy microlocal \eqref{a3} est vraie dans $\partial D$ mais pas dans $D$. Alors, $P$ a au moins une r\'esonance dans $D$ pour $h$ assez petit.
\end{corollaire}

En effet, supposons que $P$ n'a pas de r\'esonance dans $D$. Ainsi, $\chi ( P - z )^{- 1} \chi$ est holomorphe dans $D$ et continue dans $\overline{D}$. Comme on a unicit\'e du probl\`eme de Cauchy microlocal \eqref{a3}, l'implication \eqref{a5}$\Rightarrow$\eqref{a4} du th\'eor\`eme \ref{a1} montre que cette r\'esolvante tronqu\'ee est polynomialement born\'ee dans $\partial D$. Par le principe du maximum, il en est de m\^eme dans $\overline{D}$. Par l'implication \eqref{a4}$\Rightarrow$\eqref{a5} du th\'eor\`eme \ref{a1}, on en d\'eduit l'unicit\'e du probl\`eme de Cauchy microlocal \eqref{a3} dans $D$, ce qui contredit les hypoth\`eses.

On donne maintenant quelques exemples d'application du th\'eor\`eme \ref{a1}. Dans le cas non-captif, c'est-\`a-dire $K ( E_{0} )=\emptyset$, on a $\Gamma_{-} ( E_{0} ) = 
\Gamma_{+} ( E_{0} )=\emptyset$ et l'unicit\'e du probl\`eme de Cauchy microlocal est automatiquement v\'erifi\'e. Le th\'eor\`eme \ref{a1} implique donc qu'il n'y a pas de r\'esonance proche de $E_{0}$ v\'erifiant $\vert \im z \vert \leq C h \vert \ln h \vert$. On retrouve ainsi un r\'esultat de Martinez \cite{Ma02_01}. Par ailleurs, dans le cas d'un ``puits dans l'isle'', on a $\Gamma_{-} ( E_{0} ) = \Gamma_{+} ( E_{0} )=K ( E_{0} )$. Par cons\'equent, \eqref{a3} devient $( P - z ) u = v$ microlocalement pr\`es du puits (il n'y a plus de condition initiale). Soit $Q$ la r\'ealisation de Dirichlet de l'op\'erateur $P$ restreint \`a un voisinage du puits. En utilisant la r\'egularit\'e elliptique hors du puits (i.e. dans la zone physiquement interdite), on montre facilement que l'unicit\'e pour cette \'equation microlocale est vraie pour les $z$ suffisamment loins des valeurs propres de $Q$, mais est fausse sinon. Le corollaire \ref{a16} permet de conclure que les r\'esonances de $P$ co\"{i}ncident avec les valeurs propres de $Q$ modulo $\CO ( h^{\infty} )$. On retrouve ainsi une version $C^\infty$ des r\'esultats obtenus par Helffer et Sj\"{o}strand \cite{HeSj86_01} dans cette situation g\'eom\'etrique (voir aussi \cite{FuLaMa11_01,Na89_01,Na89_02}). Le th\'eor\`eme \ref{a1} peut aussi \^etre adapt\'e au cadre de l'\'etude des valeurs propres d'un op\'erateur (\ref{a0}), avec $E_0<0$, par exemple. Comme dans le cas du puits dans l'isle, la condition initiale sur $\Gamma_{-} ( E_{0} ) \setminus K ( E_{0} )$ dans (\ref{a3}) disparait. 
Le th\'eor\`eme \ref{a1} permet de caract\'eriser l'absence de valeurs propres pr\`es de $E_{0}$ par : $( P - z ) u = 0$ microlocalement pr\`es de $K ( E_{0} )$ $\Rightarrow$ $u = 0$ microlocalement pr\`es de  $K ( E_{0} )$.

Dans \cite{BoFuRaZe16_02}  et \cite{BoFuRaZe16_01}, on a donn\'e l'asymptotique des r\'esonances dans le cas d'un point fixe hyperbolique et des trajectoires homoclines/h\'et\'eroclines. Pour ce faire, on a d'abord obtenu des zones sans r\'esonances gr\^ace \`a l'implication \eqref{a5}$\Rightarrow$\eqref{a4} du th\'eor\`eme \ref{a1}. Ensuite, on a d\'emontr\'e la pr\'esence des r\'esonances en utilisant l'existence du probl\`eme de Cauchy microlocal pr\`es de l'ensemble capt\'e et des ``fonctions test'' bien choisies. Pour prouver l'unicit\'e du probl\`eme de Cauchy microlocal, on a suivi l'approche suivante : on consid\`ere une solution $u$ de \eqref{a3} avec $v = u_{0} = 0$. D'abord, on montre que $u$ poss\`ede une structure particuli\`ere (par exemple, $u$ est une distribution lagrangienne de vari\'et\'e associ\'ee bien d\'etermin\'ee). Ensuite, on projette l'\'equation \eqref{a3} sur cette structure pour obtenir une \'equation r\'eduite (par exemple, une \'equation de transport sur le symbole d'une distribution lagrangienne). Finalement, on r\'esout cette \'equation simplifi\'ee et on conclut que $u$ est nul microlocalement pr\`es de $K ( E_{0} )$.

La propagation des singularit\'es est une question classique de la th\'eorie des \'equations aux d\'eriv\'ees partielles lin\'eaires. Dans le cas pr\'esent, le d\'efaut de propagation des singularit\'es est du \`a la ``g\'eom\'etrie globale'' de l'ensemble capt\'e : le symbole $p$ peut tr\`es bien \^etre r\'eel de type principal (et donc v\'erifier le th\'eor\`eme classique d'H\"{o}rmander en tout point de $K ( E_{0} )$) bien qu'il n'y ait pas unicit\'e du probl\`eme de Cauchy microlocal \eqref{a3}. Par contre, la plupart des r\'esultats obtenus dans ce domaine concernent des ``situations irr\'eguli\`eres'' (singularit\'es des coefficients, obstacles, croisements de modes, \ldots). En adaptant la strat\'egie pr\'esent\'ee dans cette note, il devrait \^etre possible d'appliquer certains de ces r\'esultats en th\'eorie des r\'esonances.

\section{Preuve du th\'eor\`eme \ref{a1}}

L'implication \eqref{a5}$\Rightarrow$\eqref{a4} a \'et\'e d\'emontr\'ee dans la section 8 de \cite{BoFuRaZe16_01}. Dans cet article, on supposait $\vert \im z \vert \lesssim h$, mais la preuve s'\'etend sans modification au cas $\vert \im z \vert \lesssim h \vert \ln h \vert$.

On montre maintenant l'implication \eqref{a4}$\Rightarrow$\eqref{a5} par l'absurde. Supposons que \eqref{a4} soit vrai et qu'il existe une solution $u$ de
\begin{equation} \label{a6}
\left\{ \begin{aligned}
&( P - z ) u = 0 &&\text{ microlocalement pr\`es de } K ( E_{0} ) ,  \\
&u = 0 &&\text{ microlocalement pr\`es de tout point de } \Gamma_{-} ( E_{0} ) \setminus K ( E_{0} ) ,
\end{aligned} \right.
\end{equation}
avec $z \in \Omega$, $\Vert u \Vert = 1$ et telle que $u$ n'est pas nulle microlocalement pr\`es de $K ( E_{0} )$. Notons que cette assertion n'a lieu que pour une suite de valeurs de $h$ qui tend vers $0$. En particulier, il existe $\one_{K ( E_{0} )} \prec \psi \in C^{\infty}_{0} ( T^{*} \R^{n} )$ telle que
\begin{equation} \label{a7}
\Op ( \psi ) ( P - z ) u = \CO ( h^{\infty} ).
\end{equation}
La notation $f \prec g$ signifie que $g = 1$ au voisinage du support de $f$. On d\'efinit alors
\begin{equation}
v = \Op ( \varphi ) u ,
\end{equation}
pour $\varphi \in C^{\infty}_{0} ( T^{*} \R^{n} )$ avec $\one_{K ( E_{0} )} \prec \varphi \prec \psi$.

Soit $P_{\theta}$ l'op\'erateur distordu d'un angle $\theta = M h \vert \ln h \vert$ \`a partir de $P$, avec $M > 0$ assez grand. On suppose que la distortion a lieu en dehors d'un voisinage de $K ( E_{0} )$. On cherche une param\'etrixe $U ( t )$ de $e^{- i t ( P_{\theta} - z ) / h} \Op ( \psi )$ m\^eme si ce dernier op\'erateur peut ne pas \^etre bien d\'efini, $P_{\theta}$ n'\'etant pas autoadjoint. En utilisant la m\'ethode BKW, on peut construire une famille lisse d'op\'erateurs born\'es $U ( t )$ de $\CS^{\prime} ( \R^{n} )$ dans $\CS ( \R^{n} )$ telle que $U ( 0 ) = \Op ( \psi ) + \CO ( h^{\infty} )$ et
\begin{equation} \label{a15}
- i h \partial_{t} U + ( P_{\theta} - z ) U = \CO ( h^{\infty} ) ,
\end{equation}
localement uniform\'ement en $t \in \R$. De plus, $U ( t )$ est un op\'erateur int\'egral de Fourier semiclassique de transformation canonique $\exp ( t H_{p} )$ et d'ordre $h^{- C - C t}$ avec $C > 0$. En effet, comme la partie imaginaire de $P_{\theta} - z$ est d'ordre $h \vert \ln h \vert$, l'\'equation eikonale pour \eqref{a15} est la m\^eme que dans le cas autoadjoint et les \'equations de transport sont perturb\'ees par des termes d'ordre $\vert \ln h \vert$ ce qui entraine l'estimation non uniforme par rapport \`a $t$ du symbole. Comme d'habitude, la construction de $U ( t )$ se fait sur un petit intervalle de temps puis, pour tous temps, par recollement. Par ailleurs,
\begin{equation} \label{a8}
( P_{\theta} - z ) U ( t ) w = U ( t ) ( P_{\theta} - z ) w + \CO ( h^{\infty} ) ,
\end{equation}
pour toute fonction $w$ microsupport\'ee \`a l'int\'erieur de $\psi^{- 1} ( 1 )$. En effet, les deux membres de cette \'egalit\'e v\'erifient \eqref{a15} avec m\^eme donn\'ee initiale. Finalement, pour toute fonction \`a support compact $w$, microsupport\'ee dans $p^{- 1} ( E_{0} )$ mais hors de $\Gamma_{-} ( E_{0} )$, il existe $C , \varepsilon > 0$ tels que
\begin{equation} \label{a9}
U ( t ) w = \CO ( h^{- C + \varepsilon t} ) .
\end{equation}
C'est une cons\'equence du caract\`ere dissipatif de $P_{\theta} - z$ en dehors d'un compact due \`a la distortion. Plus prosa\"\i quement, cette estimation s'obtient par un argument de type Gronwall sur $\Vert U ( t ) w \Vert^{2}$.

Pour $t > 0$, on d\'efinit
\begin{equation*}
v_{t} = U ( t ) v = U ( t ) \Op ( \varphi ) u .
\end{equation*}
On montre d'abord que $v_{t}$ est un quasimode pour $P_{\theta} - z$. En utilisant \eqref{a7} et \eqref{a8}, on peut \'ecrire
\begin{align}
( P_{\theta} - z ) v_{t} &= ( P_{\theta} - z ) U ( t ) \Op ( \varphi ) u  \nonumber \\
&= U ( t ) ( P_{\theta} - z ) \Op ( \varphi ) u + \CO ( h^{\infty} ) \nonumber \\
&= U ( t ) ( P - z ) \Op ( \varphi ) u + \CO ( h^{\infty} ) \nonumber \\
&= U ( t ) \big( \Op ( \varphi ) ( P - z ) u + [ P , \Op ( \varphi ) ] u \big) + \CO ( h^{\infty} ) \nonumber \\
&= U ( t ) [ P , \Op ( \varphi ) ] u + \CO ( h^{\infty} ) . \label{a11}
\end{align}
Par r\'egularit\'e elliptique et \eqref{a7}, $[ P , \Op ( \varphi ) ] u$ est microlocalis\'e dans la surface d'\'energie $p^{- 1} ( E_{0} )$. D'autre part, \eqref{a6} et $\one_{K ( E_{0} )} \prec \varphi$ impliquent que cette fonction est nulle microlocalement pr\`es de tout point de $\Gamma_{-} ( E_{0} )$. Donc, en combinant \eqref{a9} et \eqref{a11}, il vient
\begin{equation} \label{a10}
( P_{\theta} - z ) v_{t} = \CO ( h^{- C + \varepsilon t} ) ,
\end{equation}
localement uniform\'ement en $t \in \R$.

On minore maintenant la norme de $v_{t}$. Soit $\rho$ un point du microsupport de $[ P , \Op ( \varphi ) ] u$. On vient de d\'emontrer que $p ( \rho ) = E_{0}$ et que $\rho \notin \Gamma_{-} ( E_{0} )$. Ainsi, $\exp ( t H_{p} ) ( \rho )$ s'\'echappe \`a l'infini pour $t \to + \infty$ et ne rencontre pas $K ( E_{0} )$. Par compacit\'e, il existe une fonction de troncature lisse $\chi$ ind\'ependante de $\rho$ avec $\one_{K ( E_{0} )} \prec \chi \prec\varphi$ dont le support \'evite $\exp ( [ 0 , + \infty [ H_{p} ) ( \rho )$. En particulier,
\begin{equation} \label{a12}
\Op ( \chi ) U ( t ) [ P , \Op ( \varphi ) ] u = \CO ( h^{\infty} ) ,
\end{equation}
puisque $U ( t )$ est un op\'erateur int\'egral de Fourier de transformation canonique $\exp ( t H_{p} )$. Les relations \eqref{a15}, \eqref{a8}, \eqref{a11} et \eqref{a12} impliquent
\begin{align*}
\partial_{t} \Op ( \chi ) v_{t} &= \partial_{t} \Op ( \chi ) U ( t ) \Op ( \varphi ) u  \\
&= - i h^{- 1} \Op ( \chi ) ( P_{\theta} - z ) U ( t ) \Op ( \varphi ) u + \CO ( h^{\infty} )  \\
&= - i h^{- 1} \Op ( \chi ) U ( t ) [ P , \Op ( \varphi ) ] u + \CO ( h^{\infty} )  \\
&= \CO ( h^{\infty} ) ,
\end{align*}
localement uniform\'ement en temps. En int\'egrant en temps, il vient
\begin{equation*}
\Op ( \chi ) v_{t} = \Op ( \chi ) v + \CO ( h^{\infty} ) = \Op ( \chi ) u + \CO ( h^{\infty} ) .
\end{equation*}
D'autre part, comme $u$ n'est pas nul microlocalement pr\`es de $K ( E_{0} )$, il existe $C > 0$ tel que $\Vert \Op ( \chi ) u \Vert \geq h^{C}$ pour $h$ assez petit. Ainsi, on a 
\begin{equation} \label{a13}
\Vert v_{t} \Vert \gtrsim \Vert \Op ( \chi ) v_{t} \big\Vert \geq h^{C + 1} ,
\end{equation}
pour tout $t$ fix\'e et $h$ assez petit.

D'apr\`es \eqref{a4} et la proposition D.1 de \cite{BoFuRaZe16_01}, il existe $N > 0$ tel que
\begin{equation} \label{a14}
\big\Vert ( P_{\theta} - z )^{- 1} \big\Vert \lesssim h^{- N} .
\end{equation}
En combinant \eqref{a10}, \eqref{a13} et \eqref{a14}, on obtient
\begin{equation*}
h^{C + 1} \lesssim \Vert v_{t} \Vert \leq \big\Vert ( P_{\theta} - z )^{- 1} \big\Vert \CO ( h^{ - C + \varepsilon t} ) \lesssim h^{- N - C + \varepsilon t} ,
\end{equation*}
pour tout $t$ fix\'e et $h$ assez petit. En choisissant $t$ assez grand, on obtient une contradiction en prenant $h \to 0$. Ceci conclut la preuve du th\'eor\`eme \ref{a1}.

\bibliographystyle{smfplain}
\providecommand{\bysame}{\leavevmode ---\ }
\providecommand{\og}{``}
\providecommand{\fg}{''}
\providecommand{\smfandname}{et}
\providecommand{\smfedsname}{\'eds.}
\providecommand{\smfedname}{\'ed.}
\providecommand{\smfmastersthesisname}{M\'emoire}
\providecommand{\smfphdthesisname}{Th\`ese}

\end{document}